\documentclass[12pt]{amsart}

\usepackage{mathrsfs}

\usepackage{amsmath,amsthm,amssymb}
\font\teneufm=eufm10
\font\seveneufm=eufm7
\font\fiveeufm=eufm5
\newfam\frakturfam

\textfont\frakturfam=\teneufm
\scriptfont\frakturfam=\seveneufm
\scriptscriptfont\frakturfam=\fiveeufm

\newtheorem{lm}{Lemma}
\newtheorem{theor}{Theorem}
\newtheorem{co}{Corollary}
\newtheorem{rem}{Remark}

\numberwithin{pr}{section}
\numberwithin{ex}{section}
\numberwithin{df}{section}
\numberwithin{lm}{section}
\numberwithin{theor}{section}
\numberwithin{co}{section}
\numberwithin{rem}{section}
\numberwithin{prob}{section}
\numberwithin{equation}{section}

\def\a{\alpha}
\def\b{\beta}
\def\te{\theta}
\def\g{\gamma}
\def\e{\varepsilon}
\def\dt{\delta}
\def\s{\sigma}
\def\t{\tau}

\def\Proof{{\sl Proof.}\ }

\title{Free braided nonassociative Hopf algebras and Sabinin $\tau $-algebras}

\begin{document}
\date{}
\maketitle

\begin{center}

{\bf Ualbai Umirbaev}\footnote{Wayne State University, Detroit, MI 48202, 
 USA, Eurasian National University, Astana, Kazakhstan, 
 e-mail: {\em umirbaev@math.wayne.edu}.\\  Supported by an MES grant 1226/GF3; and PAPIIT grant IN 113316, UNAM, M\'exico}  
 \\
and
\\
{\bf Vladislav Kharchenko}\footnote{FES-Cuautitl\'an, UNAM, Campo 1, M\'exico, e-mail:{\em vlad@unam.mx}.
\\
Supported by PAPIIT grant IN 113316, UNAM, and C\'atedra  PIAPI 1634, FES-C UNAM, M\'exico}
\end{center}

\begin{abstract}  Let $V$ be a linear space over a field ${\bf k}$ with a braiding  $\tau : V\otimes  V\rightarrow V\otimes  V.$
We prove that the braiding $\tau$ has a unique extension on the free nonassociative algebra ${\bf k}\{V\}$ 
freely generated by $V$ so that ${\bf k}\{V\}$ is a braided algebra. Moreover, we prove that the free 
braided algebra ${\bf k}\{V\}$ has a natural structure of a braided nonassociative Hopf algebra 
such that every element of the space of generators $V$ is primitive.
In the case of involutive braidings, $\tau^2={\rm id}$, we describe braided analogues of Shestakov-Umirbaev operations and prove that these operations are primitive operations.
We introduce a braided version of Sabinin algebras and prove that the set of all primitive elements of
a nonassociative $\tau$-algebra is a Sabinin $\tau$-algebra.
\end{abstract}

\

\noindent


\section{Introduction}
\hspace*{\parindent}

Lie theory for nonassociative products appeared as a subject of its own in the works of 
Malcev \cite{Mal55} who constructed the tangent structures corresponding to Moufang loops.

Many of the well-known generalizations of Lie algebras involve only one or two operations: Malcev algebras
have one binary bracket, Lie triple systems have one ternary bracket, Bol and Lie-Yamaguti algebras have 
one binary
and one ternary bracket, Akivis algebras have two operations --- an antisymmetric
binary bracket and a ternary bracket (related to commutator and associator) --- with only 
one identity that relates the two operations and generalizes the Jacobi identity.

For some time Akivis algebras were considered as a possible analogue of Lie algebras for nonassociative product. 
Nevertheless recent results of Shestakov and Umirbaev        \cite{SU1} demonstrate that there is an infinite set 
of independent operations, primitive polynomials of a free nonassociative algebra, that must be taken into consideration.

An important advance in the Lie theory of nonassociative products was the introduction
of a {\it hiperalgebra} by Mikheev and Sabinin, now called a {\it Sabinin algebra},
which is the most general form of the tangent structure for loops, see      \cite{MS90, S99, S00}. Lie,
Malcev, Bol, Lie-Yamaguti algebras, and Lie triple systems are specific instances
of Sabinin algebras. Sabinin algebras have an infinite set of independent operations.
There are three different natural constructions of operations in a Sabinin
algebra. Two of those constructions were devised by Sabinin and Mikheev. The third set of
operations is precisely the Shestakov-Umirbaev operations.

A construction of universal enveloping algebras, which is very similar in their properties to
usual cocommutative Hopf algebras, 
can be carried out for Bol algebras      \cite{Perez05} and, more generally, for all Sabinin
algebras      \cite{Perez07}. The role of the nonassociative Hopf algebras in the fundamental questions
of Lie theory such as integration was clarified in      \cite{MP10, MPS14}, where the authors describe the current understanding of the subject in view of the
recent works, many of which use nonassociative Hopf algebras as the main tool. 

The aim of the present paper is a developing of the quantum aspects of 
the nonassociative Lie theory, such as the
structure of primitive nonassociative polynomials whose variables form a braided space.
To this end, one should consider a free nonassociative algebra as a braided 
(nonassociative) Hopf algebra. Recall that, in his beautiful paper     \cite{Perez07}, J.M. P\'erez-Izquierdo  
elaborated a nonassociative analogue of the 
Hopf algebra concept, he called this an $H$-{\it bialgebra}. 
Instead of the antipode he had to consider two additional 
bilinear operations, the {\it left division}  and the {\it right division}:

 $$
 \setminus: H\times H\rightarrow H, \ \ \ \ \ \  / : H\times H\rightarrow H,
 $$
 $$
   (x,y)\mapsto x\setminus y, \ \ \ \ \ \ \ (x,y)\mapsto x/y,
 $$
that satisfy the following ``antipode identities"
$$
 \sum x_{(1)}\setminus (x_{(2)}y)= \e(x)y=\sum  x_{(1)} (x_{(2)}\setminus y), 
$$
$$
(y\sum x_{(1)})/x_{(2)}= \e(x)y=\sum  (y/x_{(1)}) x_{(2)}.
$$

In the case of associative Hopf algebras, the divisions reduce to $x\setminus y=S(x)y $ and $x/y=xS(y),$
where $S$ is the antipode.

The method we used in this paper is based on the analysis of the so-called {\it local action}
\cite{Gurevich} of the braid monoids similar to that given in      \cite[Section 6.2]{Khar15} for the construction 
of the free braided associative algebra. In particular, if the ground braiding is {\it involutive} 
($\tau^2={\rm id}$), then the local action of the braid monoid on
the $n$-fold tensor product $V^{\otimes n}$ reduces to the action of the symmetric
group. This allows us to use the following general {\it principle}, originally appeared
in  \cite[p. 210]{Kh05}  for generalized Gurevich Lie algebras:

\

{\it If a theorem is valid for ordinary free nonassociative  algebras and its statement in a field of braided 
algebras may be interpreted as a property of the group algebra ${\bf k}[S_n]$ under the local action, 
then this theorem is valid for free braided nonassociative  algebras as well.}

\

Indeed, if the ground braiding is involutive, then the local action of the braid monoid $B_n$ reduces to
the action of the symmetric group $S_n.$ The interpretation required in the {\it principle} reduces the 
theorem to a system of equations $F_{\lambda }(\tau_1, \tau_2, \ldots , \tau_{n-1})=0,$
where $F_{\lambda }$ are (associative) polynomials in $\tau_1, \tau_2, \ldots , \tau_{n-1}.$
Since the theorem is valid for ordinary free nonassociative algebra, it follows that 
relations $F_{\lambda }(\theta_1, \theta_2, \ldots , \theta_{n-1})=0$ are valid for the ordinary 
flip $\theta :x\otimes y\rightarrow y\otimes x.$
This implies that $F_{\lambda }(t_1, t_2, \ldots , t_{n-1})=0$ are relations in the group algebra 
$k[S_n].$ As the local action of $B_n$ reduces to
the action of the symmetric group $S_n,$ we obtain the system of required equalities
$F_{\lambda }(\tau_1, \tau_2, \ldots , \tau_{n-1})=0.$

For example, in this way Shestakov-Umirbaev theorems on primitive operations can be generalized 
 to the braided nonassociative algebra just by suitable interpretation of their statements.
 We stress, that, nevertheless, the same theorem may have a number of interpretations 
 in terms of relations, consequently, in general, a theorem may have a number of generalizations, 
 not equivalent each other in the field of braided algebra.

In Section 2, we accumulate necessary statements
on the calculations inside of the braid monoid. Then, in the third and fourth sections, we interpret  
the desired properties of a braided nonassociative algebra in terms of relations 
in braid monoids under the local action. In Section 5, we prove that there exists a unique 
braided coassociative Hopf algebra structure  on the free nonassociative algebra induced by the braiding 
of the generators space.

In Section 6, we consider in more detail symmetries (involutive braidings).
We propose a braided version for the Shestakov-Umirbaev operations and show that they remain primitive
polynomials of the braided free nonassociative Hopf algebras ($H$-bialgebra in sense of P. Izquierdo).
This allows us to define a symmetric braided version of the Sabinin algebras. 
 
Recall that in         \cite{Tak00}, M. Takeuchi proposed
a noncategorical framework in which braided bialgebras are formulated as algebras
and coalgebras with a Yang-Baxter operator (alongside compatibility conditions).
This approach is actually more general and convenient than the other approaches 
(the ground braided space need not be embedded in a braided monoidal category),
that is why we apply it here. Nevertheless, an accurate realization of the construction of nonassociative 
Hopf algebras and (in case of symmety) Sabinin algebras in terms of braided 
monoidal categories remains an interesting task. All over the paper we suppose the ground field {\bf k}
to be of characteristic zero.

\section{Braid relations}
\hspace*{\parindent}

The braid monoid $B_n$ is an associative monoid generated by braids $s_1,s_2,\ldots,s_{n-1}$ subject to the relations
\begin{equation} \label{f2.1}
s_ks_{k+1}s_k=s_{k+1}s_ks_{k+1}, \ \ s_is_j=s_js_i, \ \ \ \ 1\leq k<n-1, \ |i-j|>1.
\end{equation}
Put $[k;k]=1$ and
\begin{equation} \label{f2.2}
 \ \ [m;k]=s_{k-1}s_{k-2}\cdots s_{m+1}s_m, \ \ [k;m]=s_ms_{m+1}\cdots s_{k-2}s_{k-1}, \ \ m<k.
\end{equation}
The following relations can be found in \cite[p. 28]{Khar15}:
\begin{equation} \label{f2.3} [m;k][r;s]=[r;s][m;k],  \ \  r\leq s <m\leq k,\end{equation}  \begin{equation}
\label{f2.4} [m;k][r;m]=[r;k], \ \ r\leq m \leq k, \end{equation}  \begin{equation}
\label{f2.5} [m-1;r-1][k;t]=[k;t][m;r],  \ \  k\leq m \leq r\leq t.
\end{equation}

For any  $1\leq k\leq r<n$ put
\begin{equation} \label{f2.6}
\nu_r^{k,n}=[k;r+1] [k+1;r+2] \cdots [k+n-r-1;n].
\end{equation}

It is proved in \cite{Khar15} that
\begin{equation} \label{f2.7} \nu _r^{k,n}=[n;r][n-1;r-1]\cdots [n-r+k;k].
\end{equation}
Equalities (\ref{f2.6}) and (\ref{f2.7}) immediately imply that
\begin{equation} \label{f2.8}
\nu_{r}^{t,m}\nu_m^{m-r+t,n}= \nu_r^{t,n}=\nu_{r}^{s+1,n}\nu_s^{t,n-r+s}.
\end{equation}

Let $S_n$ be the symmetric group on the set of symbols  $\{1,2,\ldots,n\}$. A permutation $\pi\in S_n$ is called an {\em $r$-shuffle} if
$$
\pi(1)<\pi(2)<\ldots<\pi(r), \ \ \pi(r+1)<\pi(r+2)<\ldots<\pi(n).
$$
Denote by $\mathrm{Sh}_n^r$ the set of all $r$-shuffles in $S_n$.

If $\s\in \mathrm{Sh}_s^r$ and $\dt\in \mathrm{Sh}_t^p$, then we define $\s*\dt\in \mathrm{Sh}_{s+t}^{r+p}$ by
\begin{equation} \label{f2.9}
(\s*\dt)(i)=\s(i) , \ 1\leq i\leq r, \ \ \ (\s*\dt)(r+i)=s+\dt (i), \ 1\leq i\leq p.
\end{equation}
Notice that $\s*\dt$ is well-defined since every element of $\mathrm{Sh}_n^r$ is well-defined by the values on $\{1,\ldots,r\}$.

 If $\pi \in \mathrm{Sh}_n^r$, then put
\begin{equation} \label{f2.10}
[\pi]=[1;\pi(1)][2;\pi(2)]\ldots [r;\pi(r)],
\end{equation}
where $[\pi]=1$ if $r=0$. 

Let $f=f(s_1,\ldots,s_{n-1})$ be an arbitrary element of the monoid algebra of $B_n$. Then put
\begin{equation} \label{f2.11}
f_{(t)}=f(s_{t+1},\ldots,s_{t+n-1}).
\end{equation}
Notice that $f_{(t)}$ belongs to the monoid algebra of $B_{t+n}$. 
For example, if $\pi \in \mathrm{Sh}_n^r$, then
$$
[\pi]_{(t)}=[t+1;t+\pi(1)][t+2;t+\pi(2)]\ldots [t+r;t+\pi(r)]\in B_{n+t}.
$$

\begin{lm}\label{l2.1} If $\pi\in \mathrm{Sh}_r^s,$  then
$$
[\pi] \nu_r^{1,n}=\nu_r^{1,n} [\pi]_{(n-r)}.
$$
\end{lm}
\Proof
 Consider $[i;\pi(i)]$, where $1\leq i\leq s$.
Then $1\leq i<i+1\leq \pi(i)+1\leq r+1$. Applying (\ref{f2.5}), we get
$$
[i;\pi(i)]\nu_r^{1,n}=[i;\pi(i)][1;r+1][2;r+2]\ldots [n-r;n]
$$
$$    
=[i;\pi(i)][1;r+1][2;r+2]\ldots [n-r;n][i+n-r;\pi(i)+n-r] 
$$
$$
=\nu_r^{1,n} [i+n-r;\pi(i)+n-r].
$$
Hence,
$$
[\pi] \nu_r^{1,n}=[1;\pi(1)]\cdots [s;\pi(s)] \nu_r^{1,n}
$$
$$
=\nu_r^{1,n} [1+n-r;\pi(1)+n-r]\cdots [s+n-r;\pi(s)+n-r]
=\nu_r^{1,n} [\pi]_{(n-r)}. \ \ \ \ \Box
$$

\begin{lm}\label{l2.2} Let $\s\in \mathrm{Sh}_s^r$ and $\dt\in \mathrm{Sh}_t^p$. Then
$$
[\s][\dt]_{(s)}\cdot\nu_s^{r+1,p+s}=[\s*\dt].
$$
\end{lm}
\Proof By (\ref{f2.2}), (\ref{f2.6}), and (\ref{f2.10}), we have
    $$
 [\dt]_{(s)}
  \nu_{s}^{r+1,p+s}
  =[s+1;s+\dt(1)]
 \cdots [s+p;s+\dt(p)][r+1;s+1]\cdots [r+p;s+p].
  $$
Applying  (\ref{f2.3}),    we can write
  $$
 [\dt]_{(s)}
  \nu_{s}^{r+1,p+s}
  =
  ([s+1;s+\dt(1)][r+1;s+1]) \cdots ([s+p;s+\dt(p)][r+p;s+p]).
  $$
Moreover, applying (\ref{f2.4}), we get
$$
 [\dt]_{(s)}
  \nu_{s}^{r+1,p+s}
  =[r+1;s+\dt(1)]
 \cdots [r+p;s+\dt(p)].
  $$

Consequently,
$$
[\s][\dt]_{(s)}\cdot\nu_s^{r+1,p+s}
= [1;\s(1)] \cdots [r;\s(r)] [r+1;s+\dt(1)]
 \cdots [r+p;s+\dt(p)]=[\s*\dt]
  $$
  by the definition (\ref{f2.9}) of $\s*\dt$.
   $\Box$

A linear space $V$ over a field  ${\bf k}$ is called a {\it braided space} if there is a fixed  linear
mapping $\tau :V\otimes V\rightarrow V\otimes V$ (in general not necessarily invertible)
that satisfies the {\it braid relation}:
\begin{equation} \label{f2.17}
(\tau\otimes {\rm id})({\rm id}\otimes \tau )(\tau\otimes {\rm id})
=({\rm id}\otimes \tau )(\tau\otimes {\rm id})({\rm id}\otimes \tau ).
\end{equation}
{\bf Example 1.}
If $x_1, x_2, \ldots ,x_n$ is a basis of a linear space $V$,  then for arbitrary parameters
$q_{is}\in {\bf k},$ $1\leq i,s\leq n,$ the linear mapping defined by
$$
\tau: x_i\otimes x_s\mapsto q_{is}\cdot  x_s\otimes x_i
$$
is a braiding and is called a {\it diagonal braiding}.

Let $V$ be a linear space with a braiding  $\tau : V\otimes  V\rightarrow V\otimes  V$.
Consider the linear mappings
$$
\tau _i={\rm id}^{\otimes (i-1)} \otimes \tau \otimes {\rm id}^{\otimes (n-i-1)}: 
V^{\otimes n}\rightarrow V^{\otimes n}, \ \ 1\leq i<n.
$$
Due to (\ref{f2.17}), the mappings $\tau _i$ satisfy all defining relations of the braid monoid (\ref{f2.1}) :
$$
\tau_i\tau_{i+1}\tau_i=\tau_{i+1}\tau_{i}\tau_{i+1}, \ \  1\leq i<n-1; \ \ \tau_i\tau_j=\tau_j\tau_i, \ \  |i-j|>1.
$$
Therefore, $u\cdot s_i= u\tau_i$ is a well-defined action of the braid monoid $B_n$ on $V^{\otimes n}.$
 Following D. Gurevich \cite{Gurevich} we call this action a {\it local action}. 
 
It is well-known that the symmetric group $S_n$ as a monoid is defined by the generators $t_i=(i,i+1)$ (transpositions), $1\leq i<n$, and the relations
$$
t_k^2=1, \ \ t_kt_{k+1}t_k=t_{k+1}t_kt_{k+1}, \ \ t_it_j=t_jt_i, |i-j|>1.
$$

A braiding $\tau :V\otimes V\rightarrow V\otimes V$ is called {\em involutive} if $\t^2=\mathrm{id}$. 
In this case $\tau_i^2=\mathrm{id}$ and the local action of $B_n$ on $V^{\otimes n}$
 induces a {\em local action} of $S_n$ on $V^{\otimes n}$.

The best known example of an involutive braiding is the {\em ordinary flip} $\theta $ defined by   
$\theta(x\otimes y)=y\otimes x$ for all $x,y\in V$.

\begin{rem}\label{r2.1} If a braiding $\tau$ is involutive then every relation
$$
F(\tau_1,\tau_2,\ldots,\tau_{n-1})=0,
$$
where $F$ is an (associative) polynomial, holds if it holds for the ordinary flip $\tau_i=\theta_i,$ $1\leq i<n.$
\end{rem}
In fact, in the case of the ordinary flip $\tau$ the monoid generated by $\tau_1,\tau_2,\ldots,\tau_{n-1}$ is isomorphic to $S_n$ and in the case of an arbitrary involutive braiding $\tau$ this monoid 
is a homomorphic image of $S_n$.

All definitions of this section given in the language of $s_1,\ldots,s_{n-1}$ without 
any reminders will be applied to the elements $\tau_1,\ldots,\tau_{n-1}$.

\section{Braided algebras and bialgebras}
\hspace*{\parindent}

Let $V$ and $V^{\prime }$ be spaces with braidings $\tau $ and $\tau ^{\prime }$, respectively. A
linear mapping  $\varphi :V\rightarrow V^{\prime }$ is called a {\it homomorphism of braided spaces}
 if it respects the braidings; that is,
$$
\tau (\varphi \otimes \varphi )=(\varphi \otimes \varphi )\tau ^{\prime} .
$$

Let $(a\otimes b)\tau =\sum b_i\otimes a_i$ and
$(a^{\prime }\otimes b^{\prime })\tau =\sum b^{\prime }_i\otimes a^{\prime }_i$. In this case the definition of the homomorphism
 takes the form $$\sum \varphi (b_i)\otimes \varphi (a_i)=\sum \varphi (b)_i\otimes \varphi (a)_i,$$
or,  informally,  $\varphi (a_i)=\varphi (a)_i.$

An algebra $R$ with a multiplication ${\bf m}:R\otimes R\rightarrow R$
is called a {\it braided algebra} if it is a braided space and
\begin{equation} \label{f3.1}
({\bf m}\otimes {\rm id})\tau =\tau_2\tau_1({\rm id}\otimes{\bf m}),\ \ \
({\rm id}\otimes {\bf m})\tau =\tau_1\tau_2({\bf m}\otimes {\rm id}).
\end{equation}
In these formulas, as above, we use the so-called ``exponential notation"
for actions of the operators; that is, the operators in a superposition act from the left to the right.

The identity element of an algebra is usually denoted by $1$, but sometimes 
we shall use the unit mapping as well:
$$
\eta : {\bf k}\rightarrow R, \ \ \eta (\alpha)=\alpha \cdot 1.
$$

A  {\it homomorphism of braided algebras}
is a linear mapping that is both a homomorphism of algebras and braided spaces.

A {\it coalgebra} is a triple $(C,\Delta, \e)$, where $C$ is a vector space, 
$\Delta : C\rightarrow C\otimes C$ (coproduct or comultiplication)   and  
$\e : C \rightarrow {\bf k}$ (counit) are linear mappings satisfying
\begin{equation}
\label{f3.2} \Delta(\e\otimes \mathrm{id})=\mathrm{id}=\Delta(\mathrm{id}\otimes \e).
\end{equation}

A coalgebra is {\em coassocitive} if it satisfies
\begin{equation}
\label{f3.3} \Delta(\Delta\otimes \mathrm{id})=\Delta(\mathrm{id}\otimes \Delta).
\end{equation}

A  coalgebra $(C, \Delta^b ,\varepsilon )$ is called {\it braided} if it is a braided space and
\begin{equation}
\tau (\varepsilon \otimes {\rm id})=({\rm id} \otimes \varepsilon )\tau,\ \
\tau ({\rm id} \otimes \varepsilon )=(\varepsilon \otimes {\rm id})\tau ;
\label{f3.4}
\end{equation}
\begin{equation}
\tau ({\rm id}\otimes \Delta^b )=(\Delta^b \otimes {\rm id})\tau_2\tau_1,\ \ \
\tau (\Delta^b \otimes {\rm id})= ({\rm id}\otimes \Delta^b )\tau_1\tau_2.
\label{f3.5}
\end{equation}
A {\it homomorphism of braided coalgebras}   $\varphi :V\rightarrow V^{\prime }$ is a 
homomorphism of coalgebras that respects the braidings, i.e., $\varphi$ satisfies the relations 
\begin{equation}
\Delta^b(\varphi(a))=\sum_{(a)}\varphi(a^{(1)})\otimes \varphi(a^{(2)}), \ \ \varepsilon (\varphi(a))=\varepsilon (a). 
\label{f3.6}
\end{equation}

A {\it braided bialgebra} is a braided algebra and a braided coalgebra $H$
(with the same braiding) where the coproduct is an algebra homomorphism
\begin{equation} \label{f3.7}
\Delta^b :H\rightarrow H\underline{\otimes }H.
\end{equation}
Here, $H \underline{\otimes }H$ is the ordinary tensor product of spaces  with a new multiplication
\begin{equation}
(a \underline{\otimes } b)(c \underline{\otimes }d)
=\sum_i (ac_i \underline{\otimes }b_i d), \hbox{ where } (b\otimes c)\tau =\sum_i c_i\otimes b_i.
\label{f3.8}
\end{equation}

A  {\it homomorphism of braided bialgebras} is  a homomorphism of coalgebras and braided algebras.
 
In 2007 J. P\'erez-Izquierdo \cite{Perez07} elaborated the concept of a nonassociative  Hopf algebra
 and he called it an $H$-{\it bialgebra}.  An $H$-{\it bialgebra} has two additional (left and right) division  
 operations instead of the antipode. Associative $H$-bialgebras are exactly associative Hopf algebras. 
 For this reason nonassociative $H$-bialgebras are called also {\it nonassociative Hopf algebras} \cite{MPS14}.

A {\it braided nonassociative Hopf algebra} is a braided coassociative bialgebra 
$H$ with two extra bilinear operations, the left and right divisions,
$$
 \backslash: H\times H\rightarrow H, \ \ \ \ \ \  /: H\times H\rightarrow H,
 $$
 $$
   (x,y)\mapsto x\backslash y, \ \ \ \ \ \ \ (x,y)\mapsto x/y,
$$
 such that
\begin{equation}
\label{f3.9}  \sum x_{(1)}\backslash (x_{(2)}y)= \e(x)y= \sum x_{(1)} (x_{(2)}\backslash y), \end{equation}  \begin{equation}
\label{f3.10}  (y\sum x_{(1)})/x_{(2)}= \e(x)y= \sum (y/x_{(1)}) x_{(2)}.
\end{equation}

A  {\it homomorphism of braided nonassociative Hopf algebras} is  a homomorphism $\varphi$ 
of braided bialgebras that satisfies
$$
\varphi(x\backslash y)=\varphi(x)\backslash\varphi(y), \ \ \varphi(x/y)=\varphi(x)/\varphi(y).
$$

\section{Free braided nonassociative algebra}
\hspace*{\parindent}

Let $V$ be a linear space with a braiding  $\tau : V\otimes  V\rightarrow V\otimes  V$.
We fix a linear basis  $X=\{ x_i |  i\in I\} $ for $V$. The set $X^*$ of all associative words in $X$
 is a linear basis for   the free associative algebra {\bf k}$\langle X\rangle$.  This algebra
 is isomorphic to the tensor algebra $T(V)=\bigoplus_{i=0}^{\infty } V^{\otimes i}$ of the linear space $V$
with the concatenation product. We identify the words of length $m$ in $X$ with
 the corresponding tensors from $V^{\otimes m}$.  We set
$V^{\otimes 0}={\, \bf k}\cdot 1$ and  $1\otimes v=v\otimes 1=v$.

The product in the tensor algebra $T(V)$ will be denoted by ${\bf m}'$. We have 
 $(u\otimes' v){\bf m}'=u\otimes v$, where the sign $\otimes'$ is the same tensor product $\otimes$ with one additional function that separates a  pair of tensors to which the product  ${\bf m}'$ is applied.

V. Kharchenko proved that the braiding $\tau $ has a unique extension $\tau'$ on the free  algebra ${\bf k}\langle X\rangle$ so that ${\bf k}\langle X\rangle$ is a braided algebra \cite[Chapter 6]{Khar15}.
 For any $0\leq r\leq n$ denote by $\theta _r$ the linear mapping
$$\theta _r : V^{\otimes n}\rightarrow V^{\otimes r} \otimes'V^{\otimes (n-r)}$$
defined by
$$
(z_1z_2\cdots z_n)\, \theta_r=z_1z_2\cdots z_r\otimes' z_{r+1}\cdots z_{n}, \ \ \ z_i\in X.
$$
The braiding $\tau'$ is defined in \cite{Khar15} by
\begin{equation} \label{f4.1}
(u\otimes' v)\tau^{\prime } =(u\otimes v) \nu_r^{1,n}\, \theta _{n-r}, \ \ \ u\in V^{\otimes r}, \ v\in V^{\otimes (n-r)}
\end{equation}
(if $r=0$ or $r=n,$ then this definition means that
$ (1\otimes' v)\tau^{\prime }=v\otimes' 1$ or
$(u\otimes' 1)\tau^{\prime }=1\otimes' u$, respectively).

Moreover, the algebra  ${\bf k}\langle X\rangle $ has a natural structure of a braided 
Hopf algebra  such that every element of  $X$ is primitive \cite[Theorem 6.2]{Khar15}.

Denote by ${\bf k}\{ X\} = {\bf k}\{ V\}$  the free nonassociative algebra over ${\bf k}$ 
freely generated by the set $X,$ where as above $X$ is a fixed basis of the space $V.$
The product on ${\bf k}\{ X\}$ will be denoted by
$$
{\bold m} : {\bf k}\{ X\}\otimes {\bf k}\{ X\} \rightarrow {\bf k}\{ X\}.
$$

Recall that a {\it nonassociative word} is a word where the parenthesis are arranged to
show how the multiplication applies. Sometimes it is more convenient to variate a
designation of the parenthesis, for example instead of $(xy)z$ one may write $xy\cdot z,$
whereas $((z(xy))t)v$ takes the form $\{ (z\cdot xy)t\} v.$ 
A right-normed nonassociative word,
$$
u = ((\ldots ((x_1x_2)x_3)\ldots )x_m),
$$
has a simplified notation without parenthesis,
$$ 
u = x_1x_2x_3 . . . x_m.
$$
The formal definition of a nonassociative word is given in \cite[Chapter 1, Section 1]{KBKA}.

The set of all nonassociative words in the alphabet $X$ forms a linear basis for ${\bf k}\{ X\} $. 
Every nonassociative word of length $m$ has a  unique representation 
$uR=u\cdot R$, where $u$ is an associative word of length $m$ and $R$ is an 
arrangement of parenthesis.
 Because  every associative word $u\in X^*$ of length $m$ is identified with an 
 element of $V^{\otimes m}$, we can consider an arrangement of parenthesis $R$ as a linear mapping
$$
 R : V^{\otimes m} \rightarrow {\bf k}\{ V\}.
$$
We call this function a {\em parenthesis function of $m$ arguments}.  For example, if $uR=(x_1x_2)(x_3(x_4x_5))$, then $R=(\cdot\cdot)(\cdot(\cdot\cdot))$.

 Let $P$ be a free nonassociative monoid of all nonassociative words in one free variable $y$. 
 Obviously, every arrangement of parenthesis $R$ uniquely defines a nonassociative word in $y,$ 
 and we may identify by $R$ with that word.  In this way we consider every element of $P$ as a parenthesis function. 
 Denote by $P_m$ the set of all parenthesis functions of $m$ arguments. 
 We can linearly extend the action of $R\in P_m$  on  ${\bf k}\langle V\rangle$ by $V^{\otimes n}R=0$ 
  if $n\neq m$. Moreover, we can extend the action of parenthesis defined on ${\bf k}\langle V\rangle$  
  to the action of the monoid algebra ${\bf k}P$ by linearity. Obviously,
 $$
 {\bf k}\langle V\rangle \otimes {\bf k}P \rightarrow {\bf k}\{V\} \ \ \ \ \   (a\otimes p)\mapsto a\cdot p
 $$
 is an isomorphism of linear spaces, because an associative word and a parenthesis function
 uniquely defines a nonassociative word. Below we call the monoid $P$ the {\em parenthesis 
 monoid} and ${\bf k}P$ the {\em parenthesis algebra over} ${\bf k}$. 
 Notice that ${\bf k}P$ is a free nonassociative algebra in one free variable $y$.

 A record of the form $uR$, where $u\in X^*$ and $R\in P$, usually means that $u$ is an associative word of length $m$ and $R\in P_m$. Under these conditions, we have
 $$
 (uR)(vL)=(uv)(RL), \ \ \ u,v\in X^*, R,L\in P.
 $$
 \begin{theor}\label{t4.1}
 The braiding $\tau : V\otimes  V\rightarrow V\otimes  V$ has a unique extension on the free  nonassociative algebra ${\bf k}\{ X\} $ so that ${\bf k}\{ X\} $ is a braided algebra.
  \end{theor}
 \Proof Define a linear mapping
 $$ \tau^* : {\bf k}\{ X\}\otimes   {\bf k}\{ X\}\rightarrow  {\bf k}\{ X\}\otimes   {\bf k}\{ X\}$$
 by
 $$(uR\otimes   vL)\tau^*=(u\otimes' v)\tau'(L\otimes R),$$
 where $\tau'$ is the  braiding of ${\bf k}\langle X\rangle$, $u,v\in X^*$, and $R,L\in P$ defined by (\ref{f4.1}).
 This definition is equivalent to the operator equality
   \begin{equation} \label{f4.2}
 (R\otimes L)\tau^*=\tau'(L\otimes R).
   \end{equation}
Using this equality several times, we obtain
 $$
 (R\otimes L\otimes M)\tau^*_1 \tau^*_2 \tau^*_1=\tau'_1  (L\otimes R\otimes M) \tau^*_2 \tau^*_1
 =  \tau'_1 \tau'_2 (L\otimes M\otimes R) \tau^*_1
 =\tau'_1 \tau'_2 \tau'_1(M\otimes L\otimes R).
 $$
 Similarly,
 $$
 (R\otimes L\otimes M)\tau^*_2 \tau^*_1 \tau^*_2=\tau'_2 \tau'_1 \tau'_2(M\otimes L\otimes R).
 $$

 The braid relation (\ref{f2.17}) holds for $\tau',$ see \cite[Theorem 6.2]  {Khar15}.  Consequently,
  $$
 (R\otimes L\otimes M)\tau^*_1 \tau^*_2 \tau^*_1=
  (R\otimes L\otimes M)\tau^*_2 \tau^*_1 \tau^*_2
   $$
 and
    $$
 \tau^*_1 \tau^*_2 \tau^*_1=
 \tau^*_2 \tau^*_1 \tau^*_2;
   $$
 that is,  $\tau^*$ is a braiding.

Recall that  ${\bf m}'$ denotes the product in the associative algebra 
${\bf k}\langle X\rangle$ and ${\bf m}$ denotes the product in ${\bf k}\{ X\}$. 
   Applying (\ref{f4.2}), we obtain
 $$
  (uR\otimes vL\otimes wM)({\bf m}\otimes {\rm id})\tau^* = 
  (uR\cdot vL\otimes wM)\tau^*
  = (uv\otimes w)(RL\otimes M)\tau^*
  $$
  $$
  =(uv\otimes w)\tau'  (M\otimes RL)
  =(u\otimes v\otimes w)({\bf m}'\otimes {\rm id})\tau'  (M\otimes RL). 
 $$
  Using that ${\bf m}'$  and $\tau'$ satisfy (\ref{f3.1}), we have
  $$
  (uR\otimes vL\otimes wM)({\bf m}\otimes {\rm id})\tau^* =(u\otimes v\otimes w)\tau'_2\tau'_1({\rm id}\otimes {\bf m}')  (M\otimes RL).
  $$
Applying (\ref{f4.2}) several times, we also obtain
$$
(uR\otimes vL\otimes wM)\tau^*_2\tau^*_1({\rm id}\otimes {\bf m})
=
(u\otimes v\otimes w)(R\otimes L\otimes M)\tau^*_2\tau^*_1({\rm id}\otimes {\bf m})
$$
$$
=(u\otimes v\otimes w)\tau'_2\tau'_1(M\otimes R\otimes L)({\rm id}\otimes {\bf m})=
(u\otimes v\otimes w)\tau'_2\tau'_1(M\otimes RL).
$$
 Consequently,
 $$
  (uR\otimes vL\otimes wM)({\bf m}\otimes {\rm id})\tau^* =(uR\otimes vL\otimes wM)\tau^*_2\tau^*_1({\rm id}\otimes {\bf m})
  $$
  and
$$
  ({\bf m}\otimes {\rm id})\tau^* =\tau^*_2\tau^*_1({\rm id}\otimes {\bf m});
  $$
that is, the first of the relations (\ref{f3.1}) holds for $\tau^*$. One may check the second of the relations (\ref{f3.1}) similarly.   So,
${\bf k}\{ X\}$ is a braided $\tau^*$-algebra.

The uniqueness of $\tau^*$   easily follows from the relations (\ref{f3.1}).   In fact, for all $a,b,c\in {\bf k}\{ X\},$ we have
$$
(a\otimes b\otimes c)(\mathbf{m}\otimes \mathrm{id})\tau^* =(ab\otimes c)\tau^*,
$$
and
$$
(a\otimes b\otimes c)\tau^*_2\tau^*_1(\mathrm{id}\otimes\mathbf{m})
=(a\otimes (b\otimes c)\tau^*)\tau^*_1(\mathrm{id}\otimes\mathbf{m})
=(\sum_i (a\otimes c_i)\tau^*\otimes b_i )(\mathrm{id}\otimes\mathbf{m})=
\sum _{i,j} c_{ij}\otimes a_{ij}b_i,
$$
where $(b\otimes c)\tau^*=\sum_i c_i\otimes b_i$ and $(a\otimes c_i)\tau^* = \sum _jc_{ij}\otimes a_{ij}$.
Consequently, we may write the first of the relations (\ref{f3.1}) in the form
$$
(ab\otimes c)\tau^*=\sum _{i,j} c_{ij}\otimes a_{ij}b_i.
$$
This implies the uniqueness of $(ab\otimes c)\tau^*$ modulo    $(b\otimes c)\tau^*$ and $(a\otimes c_i)\tau^*$. The second of the relations (\ref{f3.1}) implies the uniqueness of  $(a\otimes bc)\tau^*$ in the same way. $\Box$

\section{Braided Hopf algebra structure on ${\bf k}\{ X\}$}

\hspace*{\parindent}

The equalities 
\begin{equation}\label{f5.1}
\Delta^b(x)=x\underline{\otimes } 1+1\underline{\otimes } x,    \ \ \ x\in X,
\end{equation}
uniquely define a homomorphism
$$
\Delta^b :{\bf k}\{ X\}\rightarrow {\bf k}\{ X\}\underline{\otimes }{\bf k}\{ X\}
$$
because ${\bf k}\{ X\}$ is freely generated by $X$. We are going to convert ${\bf k}\{ X\}$ into 
a nonassociative  Hopf algebra. 

Since $\e : {\bf k}\{ X\}\rightarrow {\bf k}$ is a homomorphism 
of algebras, it follows that (\ref{f3.2}) gives  $x \e(1)+1\e(x)=x$ for all $x\in X$. 
Consequently, $\e(x)=0$ because $\e(1)=1$. Thus, $\Delta^b$ uniquely defines the counit
\begin{equation} \label{f5.2}
\e : {\bf k}\{ X\}\rightarrow {\bf k}, \ \ \ \ \ x\mapsto 0, \ x\in X.
\end{equation}

Let $R\in P_n$ be an arbitrary $n$-ary parenthesis function and let
 $u=x_1x_2\ldots x_n$ be an arbitrary associative word of length $n$.  Let $\a=\{\a_1<\ldots<\a_r\}$ be a subset of $\{1,2,\ldots,n\}$. Denote by  ${}^{\a}(uR)$ the nonassociative word obtained from $uR$
by substituting $1$ instead of all $x_i,$ $i\neq \a_j,$ $1\leq j\leq r$. Also
${}^{\a}R$ denotes the parenthesis function corresponding to the arrangement of parenthesis in ${}^{\a}(uR)$.

 If $\pi\in \mathrm{Sh}_n^r$, then we set
 $$
{}^{\pi}(uR)={}^{\{\pi(1)<\ldots<\pi(r)\}}(uR), \ \
 (uR)^{\pi}={}^{\{\pi(r+1)<\ldots<\pi(n)\}}(uR),
$$
 and
$$
{}^{\pi}R={}^{\{\pi(1)<\ldots<\pi(r)\}}R, \ \
 R^{\pi}={}^{\{\pi(r+1)<\ldots<\pi(n)\}}R.
$$

For example, if $uR=(x_1x_2)((x_3x_4)(x_5x_6)),$ then $R=(\cdot \cdot)((\cdot \cdot)(\cdot \cdot))$. If
$$
\pi=
\left( \begin{array}{cccccc}
1 & 2 & 3 & 4 & 5 & 6 \\
2 & 6 & 1 & 3 & 4 & 5\end{array} \right)\in \mathrm{Sh}_6^2,
$$
then ${}^{\pi}(uR)=x_2x_6$, ${}^{\pi}R=(\cdot\cdot)$, $(uR)^{\pi}=(x_1((x_3x_4)x_5))$, and
$R^{\pi}=(\cdot((\cdot \cdot)\cdot))$.
Of course,  one may consider an arbitrary nonassociative word instead of $uR.$
For example, if $w=((x_6x_1)x_1)((x_3x_4)x_3),$
 then
${}^{\pi}w=x_1x_3,$ and
$w^{\pi}=(x_6x_1)(x_3x_4)$.

\begin{lm}\label{l5.1} Let $R$ be an arbitrary parenthesis function of $n$ arguments. Then
\begin{equation} \label{f5.3}
R \Delta^b= \sum _{r=0}^n \sum _{\pi\in \mathrm{Sh}_n^r} [\pi] ({}^{\pi}R  \underline{\otimes } R^ {\pi}).
\end{equation}
\end{lm}
\Proof  Put $u=x_1x_2\ldots x_n$.
If $n=1$, then $R=(\cdot)$ and (\ref{f5.3}) is obviously true. If $n=2$, then $R=(\cdot \cdot)$.  Using (\ref{f5.1}), we get
$$
u\cdot R\Delta^b=((x_1x_2))\Delta^b=x_1x_2\underline{\otimes } 1+x_1\underline{\otimes } x_2+(x_1\underline{\otimes } x_2)\tau +1\underline{\otimes } x_1x_2.
$$
The right hand-side of (\ref{f5.3}) applied to $u$ gives
$$
 u\cdot ({\rm id} (1\underline{\otimes } (\cdot \cdot ))+[1;1]((\cdot)  \underline{\otimes } (\cdot))
+[1;2]((\cdot)  \underline{\otimes } (\cdot))+[1;1][2;2]((\cdot\cdot)  \underline{\otimes } 1 ))
$$
$$
=   1\underline{\otimes } x_1x_2 +x_1\underline{\otimes } x_2
+(x_1\underline{\otimes } x_2)\tau +x_1x_2\underline{\otimes } 1.
$$
This proves (\ref{f5.3}) for $n=2$.

We prove (\ref{f5.3}) by induction on $n$. If $n\geq 2$, then $R$ has a unique  representation
in the form $R=LM,$  and, respectively, $uR=vL\cdot wM$ for some nontrivial associative words  $v$ and $w$. Suppose that $L$ and $M$ are functions of $s$ and $t$ arguments, respectively. We may assume that (\ref{f5.3}) is true for $L$ and $M$ since $s,t<n$. Consequently,
$$
u\cdot R\Delta^b=\Delta^b(uR)=\Delta^b(vL)\underline{\otimes}\Delta^b(wM)
$$
$$
= (\sum _{r=0}^s \sum _{\s\in \mathrm{Sh}_s^r} v[\s] ({}^{\s}L  \underline{\otimes} L^ {\s}))\underline{\otimes}
 (\sum _{p=0}^t \sum _{\dt\in \mathrm{Sh}_t^p} w[\dt]({}^{\dt}M  \underline{\otimes} M^ {\dt})).
$$
Suppose that
$$
v[\s] =\sum_i (v_{1i}  {\otimes } v_{2i}), \ \ \
w[\dt]
=\sum_j (w_{1j}  {\otimes } w_{2j}).
$$
In this case
$$
v[\s] ({}^{\s}L  \underline{\otimes } L^ {\s})
=\sum_i (v_{1i} ({}^{\s}L) \underline{\otimes } v_{2i}L^ {\s}), \ \ \ \
w[\dt]({}^{\dt}M  \underline{\otimes } M^ {\dt})
=\sum_j (w_{1j}({}^{\dt}M)  \underline{\otimes } w_{2j}M^ {\dt}).
$$
By the definition of the product in ${\bf k}\{ X\}\underline{\otimes }{\bf k}\{ X\}$,  we get
$$
(v_{1i} ({}^{\s}L) \underline{\otimes } v_{2i}L^ {\s})\underline{\otimes } (w_{1j}({}^{\dt}M)  \underline{\otimes } w_{2j}M^ {\dt})
=(v_{1i} ({}^{\s}L))(v_{2i}L^ {\s}\underline{\otimes } w_{1j}({}^{\dt}M))^{\tau^*}(w_{2j}M^ {\dt}),
$$
where $a(b\underline{\otimes } c)d$ equals the element $ab \underline{\otimes } cd$. Moreover, we have
$$
(v_{2i}L^ {\s}\underline{\otimes } w_{1j}({}^{\dt}M))^{\tau^*}=(v_{2i}\otimes
w_{1j})^{\tau'}({}^{\dt}M \underline{\otimes } L^ {\s})
$$
by the definition of $\tau^*$.
Therefore,
$$
(v_{1i} ({}^{\s}L) \underline{\otimes } v_{2i}L^ {\s})\underline{\otimes } (w_{1j}({}^{\dt}M)  \underline{\otimes } w_{2j}M^ {\dt})
=(v_{1i}\otimes (v_{2i}\otimes
w_{1j})^{\tau'}\otimes w_{2j}) (({}^{\s}L)({}^{\dt}M)\underline{\otimes } L^ {\s}M^ {\dt}).
$$
By the definition of $\tau'$, we have
$$
(v_{2i}\otimes
w_{1j})^{\tau'}=(v_{2i}\otimes
w_{1j})\nu_{s-r}^{1,p+s-r}\theta_p.
$$
Then
$$
v_{1i}\otimes (v_{2i}\otimes
w_{1j})^{\tau'}\otimes w_{2j}=(v_{1i}\otimes v_{2i}\otimes
w_{1j}\otimes w_{2j})\nu_{s}^{r+1,p+s}\theta_{r+p},
$$
and
$$
(v_{1i} ({}^{\s}L) \underline{\otimes } v_{2i}L^ {\s})\underline{\otimes }(w_{1j}({}^{\dt}M)  
\underline{\otimes } w_{2j}M^ {\dt})
$$
$$
=(v_{1i}\otimes v_{2i}\otimes
w_{1j}\otimes w_{2j})\nu_{s}^{r+1,p+s} (({}^{\s}L)({}^{\dt}M)\underline{\otimes } L^ {\s}M^ {\dt}).
$$
Consequently,
$$
v[\s] ({}^{\s}L  \underline{\otimes } L^ {\s})\underline{\otimes } w[\dt]({}^{\dt}M  \underline{\otimes } M^ {\dt})
= \sum_{i,j} (v_{1i} ({}^{\s}L) \underline{\otimes } v_{2i}L^ {\s})\underline{\otimes } (w_{1j}({}^{\dt}M)  \underline{\otimes } w_{2j}M^ {\dt})
$$
$$
=\sum_{i,j}(v_{1i}\otimes v_{2i}\otimes
w_{1j}\otimes w_{2j})\nu_{s}^{r+1,p+s} (({}^{\s}L)({}^{\dt}M)\underline{\otimes } L^ {\s}M^ {\dt}).
$$
Furthermore,
$$
\sum_{i,j}(v_{1i}\otimes v_{2i}\otimes
w_{1j}\otimes w_{2j})=(\sum_i v_{1i}\otimes v_{2i})\otimes (\sum_j w_{1j}\otimes w_{2j})
$$
$$
=v[\s]
\otimes w [\dt] =v[1;\s(1)] \cdots [r;\s(r)]
\otimes w [1;\dt(1)]\cdots [p;\dt(p)] 
$$
$$
=(v\otimes w)  [1;\s(1)] \cdots [r;\s(r)]  [s+1;s+\dt(1)]\cdots [s+p;s+\dt(p)]
$$
$$
=(v\otimes w)[\s] [\dt]_{(s)},
$$
and
$$
\sum_{i,j}(v_{1i}\otimes v_{2i}\otimes
w_{1j}\otimes w_{2j})\nu_{s}^{r+1,p+s}
=(v\otimes w) [\s] [\dt]_{(s)}\nu_{s}^{r+1,p+s}=(v\otimes w) [\s*\dt]
  $$
  by Lemma \ref{l2.2}. Hence,
$$
v[\s] ({}^{\s}L  \underline{\otimes } L^ {\s})\underline{\otimes } w[\dt]({}^{\dt}M  
\underline{\otimes } M^ {\dt})=(v\otimes w)  [\s*\dt] (({}^{\s}L)({}^{\dt}M)\underline{\otimes } L^ {\s}M^ {\dt}).
$$
Consequently,
$$
u\cdot R\Delta^b=\Delta^b(uR)
= \sum _{r=0}^s \sum _{p=0}^t \sum _{\s\in \mathrm{Sh}_s^r}\sum _{\dt\in \mathrm{Sh}_t^p}(v\otimes w) \cdot [\s*\dt] (({}^{\s}L)({}^{\dt}M)\underline{\otimes } L^ {\s}M^ {\dt}) .
$$
Notice that
$$
u=v\otimes w, ({}^{\s}L)({}^{\dt}M)={}^{\s*\dt}(LM), L^ {\s}M^ {\dt}=(LM)^ {\s*\dt}.
 $$
 Finally,
 $$
u\cdot R\Delta^b=\sum _{r=0}^s \sum _{p=0}^t \sum _{\s\in \mathrm{Sh}_s^r}
\sum _{\dt\in \mathrm{Sh}_t^p} u \cdot [\s*\dt] ({}^{\s*\dt}(LM) \underline{\otimes } (LM)^ {\s*\dt} ).
$$
 Consider an arbitrary shuffle $\pi \in \mathrm{Sh}_n^q$. For any given $s,t$ such that $s+t=n$,
 we can uniquely define integers $r,p$ and shuffles $\s\in \mathrm{Sh}_s^r$, $\dt\in \mathrm{Sh}_t^p$ 
 such that $\pi=\s*\dt$ and $q=r+p$. In fact, $r$ is the greatest integer such that $r\leq q$ and $\pi(r)\leq s,$ 
 whereas $\s$ is uniquely defined by $\s(i)=\pi(i),$ $1\leq i\leq r$.
 Then $p=q-r,$ and $\dt(i)=\pi(r+i),$ $1\leq i\leq p$.

 Using this decomposition of shuffles, we can replace
 $\sum _{r=0}^s \sum _{p=0}^t \sum _{\s\in \mathrm{Sh}_s^r}\sum _{\dt\in \mathrm{Sh}_t^p} $ with
  $\sum _{q=0}^n \sum _{\pi\in \mathrm{Sh}_n^q}$ and
 $\s*\dt$ with $\pi$ in the last formula for $u\cdot R\Delta^b$. This gives (\ref{f5.3}). $\Box$

\begin{lm}\label{l5.2}
The coproduct $\Delta^b$ is coassociative.
\end{lm}
\Proof We have to check that (\ref{f3.3}) holds. Let  $R\in P$ be an arbitrary 
parenthesis function of $n$ arguments. By (\ref{f5.3}), we have
$$
R \Delta^b(\Delta^b \otimes \mathrm{id})= \sum _{p=0}^n \sum _{\pi\in \mathrm{Sh}_n^p} [\pi] ({}^{\pi}R  \underline{\otimes } R^ {\pi})(\Delta^b \otimes \mathrm{id})
$$
$$
=\sum _{p=0}^n \sum _{\pi\in \mathrm{Sh}_n^p} [\pi] ({}^{\pi}R \Delta^b \underline{\otimes } R^ {\pi})=
\sum _{p=0}^n \sum _{\pi\in \mathrm{Sh}_n^p} [\pi] (\sum _{r=0}^p 
\sum _{\s\in \mathrm{Sh}_p^r} [\s] ({}^{\s}({}^{\pi}R)  \underline{\otimes } 
({}^{\pi}R)^ {\s}) \underline{\otimes } R^ {\pi})
$$
$$
=\sum _{p=0}^n \sum _{r=0}^p\sum _{\pi\in \mathrm{Sh}_n^p} 
\sum _{\s\in \mathrm{Sh}_p^r}[\pi] [\s] ({}^{\s}({}^{\pi}R)  \underline{\otimes } 
({}^{\pi}R)^ {\s} \underline{\otimes } R^ {\pi})
$$
$$
=\sum _{r+s\leq n} \sum _{\pi\in \mathrm{Sh}_n^{r+s}} 
\sum _{\s\in \mathrm{Sh}_{r+s}^r}[\pi] [\s] ({}^{\s}({}^{\pi}R)  \underline{\otimes } 
({}^{\pi}R)^ {\s} \underline{\otimes } R^ {\pi}).
$$
Similarly,
$$
R \Delta^b(\mathrm{id}\otimes \Delta^b)= \sum _{r=0}^n \sum _{\te\in \mathrm{Sh}_n^r} [\te] ({}^{\te}R  \underline{\otimes } R^ {\te})(\mathrm{id}\otimes \Delta^b)
$$
$$
=\sum _{r=0}^n \sum _{\te\in \mathrm{Sh}_n^r} [\te] ({}^{\te}R  \underline{\otimes } R^ {\te}\Delta^b)
=\sum _{r=0}^n \sum _{\te\in \mathrm{Sh}_n^r} [\te] ({}^{\te}R  \underline{\otimes } 
\sum _{s=0}^{n-r} \sum _{\dt\in \mathrm{Sh}_{n-r}^s} [\dt] ({}^{\dt}(R^ {\te})  \underline{\otimes } (R^ {\te})^ {\dt}))
$$
$$
=\sum _{r=0}^n \sum _{s=0}^{n-r}\sum _{\te\in \mathrm{Sh}_n^r} 
 \sum _{\dt\in \mathrm{Sh}_{n-r}^s}[\te] ({}^{\te}R  \underline{\otimes }  [\dt] ({}^{\dt}(R^ {\te})  
 \underline{\otimes } (R^ {\te})^ {\dt}))
$$
$$
=\sum _{r=0}^n \sum _{s=0}^{n-r}\sum _{\te\in \mathrm{Sh}_n^r}  \sum _{\dt\in \mathrm{Sh}_{n-r}^s}[\te]
[r+1;r+\dt(1)]\cdots [r+s;r+\dt(s)]
({}^{\te}R  \underline{\otimes }  {}^{\dt}(R^ {\te})  \underline{\otimes } (R^ {\te})^ {\dt})
$$
$$
=\sum _{r=0}^n \sum _{s=0}^{n-r}\sum _{\te\in \mathrm{Sh}_n^r}  \sum _{\dt\in \mathrm{Sh}_{n-r}^s}[\te]
[\dt]_{(r)}
({}^{\te}R  \underline{\otimes }  {}^{\dt}(R^ {\te})  \underline{\otimes } (R^ {\te})^ {\dt})
$$
$$
=\sum _{r+s\leq n}\sum _{\te\in \mathrm{Sh}_n^r}  \sum _{\dt\in \mathrm{Sh}_{n-r}^s}[\te]
[\dt]_{(r)}
({}^{\te}R  \underline{\otimes }  {}^{\dt}(R^ {\te})  \underline{\otimes } (R^ {\te})^ {\dt}).
$$
In order to prove the statement of the lemma, we need to prove
\begin{equation} \label{f5.4}
\sum _{r+s\leq n} \sum _{\pi\in \mathrm{Sh}_n^{r+s}} \sum _{\s\in \mathrm{Sh}_{r+s}^r}[\pi] [\s] ({}^{\s}({}^{\pi}R)  \underline{\otimes } ({}^{\pi}R)^ {\s} \underline{\otimes } R^ {\pi})
\end{equation}
$$
=\sum _{r+s\leq n}\sum _{\te\in \mathrm{Sh}_n^r}  \sum _{\dt\in \mathrm{Sh}_{n-r}^s}[\te]
[\dt]_{(r)}
({}^{\te}R  \underline{\otimes }  {}^{\dt}(R^ {\te})  \underline{\otimes } (R^ {\te})^ {\dt}).\nonumber
$$

Let $\a=\{\a_1<\ldots<\a_r\}$, $\b=\{\b_1<\ldots<\b_s\}$, and $\g=\{\g_1<\ldots<\g_t\}$ be a partition of $\{1,2,\ldots,n\}$ in three subsets. Define $\pi\in \mathrm{Sh}_n^{r+s}$ by
$$
\{\pi(1)<\ldots<\pi(r+s)\}=\a\cup\b.
$$
Obviously, $R^ {\pi}={}^{\g}R$. Let $i_1<\ldots<i_r$ be integers such that
$$
\{\pi(i_1)<\ldots<\pi(i_r)\}=\a.
$$
Define $\s\in \mathrm{Sh}_{r+s}^r$ by
$$
\{\s(1)<\ldots<\s(r)\}=\{i_1<\ldots<i_r\}.
$$
It is easy to see that ${}^{\s}({}^{\pi}R)={}^{\a}R$ and $({}^{\pi}R)^{\s}={}^{\b}R$. Notice also that
$$
(\pi\s)(i)=\a_i,  1\leq i\leq r;\ \ \ (\pi\s)(r+j)=\b_j,    1\leq j\leq s.
$$

Consider the increasing sequence $\pi(1),\ldots,\pi(r+s)$. This sequence has two increasing subsequences $\a_1,\ldots,\a_r$ and $\b_1,\ldots,\b_s$. We turn $\pi(1),\ldots,\pi(r+s)$ into $\a_1,\ldots,\a_r,\b_1,\ldots,\b_s$ by moving the elements $\b_s,\b_{s-1},\ldots,\b_1$ in the line to the right, i.e., first move $\b_s$ to the rightmost position applying transpositions, then $\b_{s-1}$ to the second rightmost position, and so on. Let $j_i$ be the number of transpositions applied to $\b_i$ during this procedure. We are going to prove that
\begin{equation} \label{f5.5}
[\pi] [\s]=[1;\a_1]\ldots[r;\a_r][r+1;\b_1+j_1]\ldots[r+s;\b_s+j_s].
\end{equation}

We have,
$$
[\pi]=[1;\pi(1)]\ldots[r+s;\pi(r+s)], \ \ \ [\s]=[1;\s(1)]\ldots[r;\s(r)].
$$
We demonstrate how to calculate $[1;\pi(1)]\ldots[r;\pi(r)] [1;\s(1)]$.  First we move $[1;\s(1)]$ to the left until it meets   $[\s(1);(\pi\s)(1)] $ using (\ref{f2.3}) and change $[\s(1);(\pi\s)(1)] [1;\s(1)]$   by  $[1;(\pi\s)(1)] $   using (\ref{f2.4}).  Then we move $[1;(\pi\s)(1)] $ to the leftmost position by (\ref{f2.5}). We obtain an equality of the form
$$
[1;\pi(1)]\ldots[r+s;\pi(r+s)] [1;\s(1)]=[1;(\pi\s)(1)][2;d_2]\ldots [r+s;d_{r+s}].
$$
We have $2\leq d_2<\ldots<d_{r+s}$ despite replacing of some $[u;v]$ by $[u+1;v+1]$
 using (\ref{f2.5}) because $[\s(1);(\pi\s)(1)]$ is disappeared in $[2;d_2]\ldots [r+s;d_{r+s}]$.
 Moreover, all multipliers  of $[1;\pi(1)]\ldots[r+s;\pi(r+s)]$ after $[\s(1);(\pi\s)(1)]$ remain unchanged
 in   $[2;d_2]\ldots [r+s;d_{r+s}]$. This allows us to calculate $[2;d_2]\ldots [r+s;d_{r+s}][2;\s(2)]$ in the same manner. Continuing these discussions we obtain
 $$
 [\pi] [\s]=[1;(\pi\s)(1)]\cdots [r;(\pi\s)(r)] [r+1;\mu_{r+1}]\ldots [r+s;\mu_{r+s}],
    $$
where $r+1\leq \mu_{r+1}<\ldots<  \mu_{r+s}\leq n$.  Notice that $[r+i;\mu_{r+i}]$ appeared instead of
$[\s(r+i);(\pi\s)(r+i)]=[\s(r+i);\b_i]$. In order to move
$[\s(r+i);\b_i]$ to the place of $[r+i;\mu_{r+i}]$ we apply relations of the form (\ref{f2.5}) exactly
$j_i$ times. Consequently, $\mu_{r+i}=\b_i+j_i$. This proves (\ref{f5.5}).

Define $\te\in \mathrm{Sh}_n^r$ by $\{\te(1)<\ldots<\te(r)\}=\a$. Then ${}^{\te}R={}^{\a}R$ and
$\{\te(r+1)<\ldots<\te(n)\}=\b\cup\g$. Let $k_1<\ldots<k_s$ be integers such that
$$
\{\te(r+k_1)<\ldots<\te(r+k_s)\}=\{\b_1<\ldots<\b_s\}.
$$
Define $\dt\in \mathrm{Sh}_{s+t}^s$ by $\{\dt(1)<\ldots<\dt(s)\}=\{k_1<\ldots<k_s\}$. It is easy to see that
${}^{\dt}(R^ {\te})={}^{\b}R$ and ${}^{\pi}(R^ {\te})={}^{\g}R$.

Notice that there are $k_i-1$ elements in $\b\cup\g$ less than $\b_i$ by the definition of $k_i$. There are
$r-j_i$ elements in $\a$ less than $\b_i$ by the definition of $j_i$. Therefore,
$$
\b_i=k_i-1+r-j_i+1=k_i+r-j_i, 1\leq i\leq s.
$$
Then,
$$
[\dt]_{(r)}=[r+1;r+\dt(1)]\cdots [r+s;r+\dt(s)] =
[r+1;r+k_1]\cdots [r+s;r+k_s]
$$
$$
=
[r+1;\b_1+j_1]\cdots [r+s;\b_s+j_s]
$$
and
$$
[\te][\dt]_{(r)}
=[1;\a_1]\cdots[r;\a_r][r+1;\b_1+j_1]\cdots [r+s;\b_s+j_s]=[\pi][\s]
$$
by (\ref{f5.5}).

Notice that $\a,\b,\g$ uniquely define 
$(\pi,\s)\in \mathrm{Sh}_n^{r+s} \times \mathrm{Sh}_{r+s}^r$ and 
$(\te,\dt)\in \mathrm{Sh}_n^r \times \mathrm{Sh}_{s+t}^s.$ 
This establishes a one-to-one correspondence between 
$\mathrm{Sh}_n^{r+s} \times \mathrm{Sh}_{r+s}^r$ and $\mathrm{Sh}_n^r \times \mathrm{Sh}_{s+t}^s$.

Consider
$$
A=[\pi] [\s] ({}^{\s}({}^{\pi}R)  \underline{\otimes } ({}^{\pi}R)^ {\s} \underline{\otimes } R^ {\pi})
$$
and
$$
A'=[\te] [\dt]_{(r)} ({}^{\te}R  \underline{\otimes }  {}^{\dt}(R^ {\te})  \underline{\otimes } (R^ {\te})^ {\dt}).
$$
We already have proved that $A=A'$. Notice that $A$ is a summand of the left hand side of (\ref{f5.4}) and
$A'$ is a summand of the right hand side of (\ref{f5.4}). The correspondence between the pairs $(\pi,\s)$ and $(\te,\dt)$ for all $r,s$ gives a one-to-one correspondence between the summands $A$ and $A'$ of the left and right hand sides of (\ref{f5.4}). This proves  that (\ref{f5.4}) is valid. $\Box$

  \begin{lm}\label{l5.3}
The braided nonassociative algebra ${\bf k}\{X\}$ with the coproduct $\Delta^b$ defined 
by $($\ref{f5.1}$)$ and the counit $\e$ define by $($\ref{f5.2}$)$  is a braided bialgebra.
\end{lm}
\Proof Let's check the first of the relations (\ref{f3.5}). Let $R\in P_r$ and $L\in P_{n-r}$. Then
$$
(R \underline{\otimes} L)\t^*(\mathrm{id}\underline{\otimes} \Delta^b)
=\nu_r^{1,n}(L \underline{\otimes} R)(\mathrm{id}\underline{\otimes} \Delta^b)
=\nu_r^{1,n}(L \underline{\otimes} R\Delta^b)
$$
$$
=\nu_r^{1,n}(L \underline{\otimes} \sum _{s=0}^r \sum _{\pi\in \mathrm{Sh}_r^s} [\pi] ({}^{\pi}R  \underline{\otimes } R^ {\pi}))=\sum _{s=0}^r \sum _{\pi\in \mathrm{Sh}_r^s}\nu_r^{1,n}[\pi]_{(n-r)}(L \underline{\otimes}   {}^{\pi}R  \underline{\otimes } R^ {\pi})
$$
and
$$
(R \underline{\otimes} L)(\Delta^b\underline{\otimes} \mathrm{id})\t_2^*\t_1^*
=(R\Delta^b \underline{\otimes} L)\t_2^*\t_1^*
=(\sum _{s=0}^r \sum _{\pi\in \mathrm{Sh}_r^s} [\pi] ({}^{\pi}R  \underline{\otimes } 
R^ {\pi}) \underline{\otimes} L)\t_2^*\t_1^*
$$
$$
=(\sum _{s=0}^r \sum _{\pi\in \mathrm{Sh}_r^s} 
[\pi] ({}^{\pi}R  \underline{\otimes } R^ {\pi}) \underline{\otimes} L)\t_2^*\t_1^*
= \sum _{s=0}^r \sum _{\pi\in \mathrm{Sh}_r^s} [\pi] ({}^{\pi}R  \underline{\otimes } 
R^ {\pi} \underline{\otimes} L)\t_2^*\t_1^*.
$$
Notice that
$$
({}^{\pi}R  \underline{\otimes } R^ {\pi} \underline{\otimes} L)\t_2^*\t_1^*=
({}^{\pi}R  \underline{\otimes } (R^ {\pi} \underline{\otimes} L)\t^*)\t_1^*
=({}^{\pi}R  \underline{\otimes } \nu_{r-s}^{1,n-s}(L \underline{\otimes} R^ {\pi}))\t_1^*
$$
$$
=\nu_{r}^{s+1,n}({}^{\pi}R  \underline{\otimes } L \underline{\otimes} R^ {\pi})\t_1^*
=\nu_{r}^{s+1,n}(({}^{\pi}R  \underline{\otimes } L)\t^* \underline{\otimes} R^ {\pi})
$$
$$
=\nu_{r}^{s+1,n}(\nu_s^{1,n-r+s}(L \underline{\otimes } {}^{\pi}R) \underline{\otimes} R^ {\pi})
=\nu_{r}^{s+1,n}\nu_s^{1,n-r+s}(L \underline{\otimes }  {}^{\pi}R \underline{\otimes} R^ {\pi})
$$
$$
=\nu_{r}^{1,n}(L \underline{\otimes }  {}^{\pi}R \underline{\otimes} R^ {\pi})
$$
by (\ref{f2.7}).
Consequently,
$$
(R \underline{\otimes} L)(\Delta^b\underline{\otimes} \mathrm{id})\t_2^*\t_1^*=
\sum _{s=0}^r \sum _{\pi\in \mathrm{Sh}_r^s} [\pi] \nu_r^{1,n}
(L \underline{\otimes }  {}^{\pi}R \underline{\otimes} R^ {\pi}).
$$

Notice that $\nu_r^{1,n}=[1;r+1][2;r+2]\ldots [n-r;n]$. Consider $[i;\pi(i)]$, where $1\leq i\leq s$.
Then $1\leq i<i+1\leq \pi(i)+1\leq r+1$. Consequently, applying Lemma \ref{l2.1}, we get
$$
[i;\pi(i)]\nu_r^{1,n}=\nu_r^{1,n} [i+n-r;\pi(i)+n-r].
$$
Hence
$$
[\pi] \nu_r^{1,n}=[1;\pi(1)]\cdots [s;\pi(s)] \nu_r^{1,n}
$$
$$
=\nu_r^{1,n} [1+n-r;\pi(1)+n-r]\cdots [s+n-r;\pi(s)+n-r]=\nu_r^{1,n} [\pi]_{(n-r)}.
$$
This proves that
$$
(R \underline{\otimes} L)\t^*(\mathrm{id}\underline{\otimes} \Delta^b)=
(R \underline{\otimes} L)(\Delta^b\underline{\otimes} \mathrm{id})\t_2^*\t_1^*;
$$
that is,
$$
\t^*(\mathrm{id}\underline{\otimes} \Delta^b)=
(\Delta^b\underline{\otimes} \mathrm{id})\t_2^*\t_1^*.
$$
Similarly one can check the second of the relations (\ref{f3.5}). The relations (\ref{f3.2}) and (\ref{f3.4}) 
hold trivially. $\Box$
\begin{theor}\label{t5.1}
There is a unique braided nonassociative Hopf algebra structure on ${\bf k}\{ X\}$ such that all elements of the generator set $X$ are primitive.
\end{theor}
\Proof It remains only to show the existence and uniqueness of the left and right divisions.  
We consider only the left one. 
We shall prove that the first equation of (\ref{f3.9}) uniquely defines $\backslash$. In fact, $1\backslash y=y$ since $\Delta^b(1)=1\otimes 1$ and $\varepsilon(1)=1$.

Let $u$ be an arbitrary nonassociative word of length $m\geq 1$. By Lemma \ref{l5.1}, $u\Delta^b$ can be uniquely written in the form
$$
u\Delta^b= u \underline{\otimes} 1 +\sum \a_i a_i \underline{\otimes} c_i,
$$
where $a_i$, $c_i$ are nonassociative words and length of each $a_i$ is less than $m$. 
Using the first equality of (\ref{f3.9}) we obtain
$$
u\backslash y=-\sum \a_i a_i \backslash (c_iy)
$$
because $\varepsilon(u)=0$.
This provides a recursive definition of $u\backslash y$ by the length of $u$.

The vector space $\mathrm{Hom}({\bf k}\{ X\}, \mathrm{End}({\bf k}\{ X\}))$ with the convolution product
$$
(\varphi*\psi)=\sum \varphi(x_{(1)})\psi(x_{(2)})
$$
is an associative algebra with the identity $\iota : x\mapsto \varepsilon(x)\mathrm{Id}$ \cite[p. 839]{Perez07}.

Consider
$$
{\mathfrak L }: {\bf k}\{ X\}\rightarrow \mathrm{End}({\bf k}\{ X\}), 
\ \ x\mapsto {\mathfrak L }_x, \ {\mathfrak L }_x(a)=xa,
$$
and
$$
 \backslash: {\bf k}\{ X\}\rightarrow \mathrm{End}({\bf k}\{ X\}), \ \ x\mapsto \backslash_x,
 \ \backslash_x(a)=x\backslash a.
$$
The first equality of (\ref{f3.9}) exactly means that $\backslash*{\mathfrak L }=\iota$;
 that is, $\backslash$ is a left inverse to ${\mathfrak L }.$

Using the second equality of (\ref{f3.9}) we can similarly define $\backslash'$ and check that
${\mathfrak L }*\backslash'=\iota$, i.e., $\backslash'$ is a right inverse to 
${\mathfrak L }$. Consequently, $\backslash=\backslash'$. $\Box$

\section{Primitive elements}

\hspace*{\parindent}

Let $B$ be an arbitrary braided (nonassociative) Hopf algebra with the braiding $\tau : B\otimes B\rightarrow B\otimes B$. Put
$$
{\bf m}_i={\rm id}^{\otimes (i-1)} \otimes {\bf m} \otimes {\rm id}^{\otimes (n-i-1)}: 
B^{\otimes n}\rightarrow B^{\otimes n-1}, \ \ 1\leq i\leq n-1,
$$
where ${\bf m} : B\otimes B\rightarrow B$ is the multiplication.

For example, if $u=u_1\otimes u_2\otimes\ldots\otimes u_k\in B^k$, then
$$
u{\bf m}_1^{k-1}=(\ldots(u_1u_2)\ldots u_k)
$$
is the right-normed product of the elements $u_1,u_2,\ldots,u_k$. In this section we denote by $R_k$ the parenthesis function such that
$$
uR_k=(\ldots(u_1u_2)\ldots u_k);
$$
that is, $R_k$ is the {\it right-normed $k$-ary parenthesis function}.

An $n$-ary operation $T$ on $B$ is called {\em primitive} if $T(u_1,\ldots,u_n)$ is primitive for any primitive elements $u_1,\ldots,u_n\in B$. For example, the commutator $[x,y]=xy-yx$ and the associator $(x,y,z)=xyz-x(yz)$ are primitive operations \cite[Lemma 8.2]{Khar15} in the case of the ordinary flip 
($\tau(x\otimes y)=(y\otimes x)$). I.P. Shestakov and U.U. Umirbaev \cite[p. 539]{SU1} introduced
a system  of primitive operations $p_{m,n}$ for all $m,n\geq 1$ in the following way.

Given $m,n\geq 1,$ let ${\bf U} = (u_1,u_2, \ldots ,u_m)$ and ${\bf V} = (v_1,v_2, \ldots ,v_n)$ 
be sequences of nonassociative polynomials, and let 
$U = u_1u_2 \cdots u_m,$ $V = v_1v_2 \cdots v_n$ be the corresponding
right-normed products. The operations are defined inductively:
$$
p({\bf U};{\bf V};w) = (U,V,w)-\sum U_{(1)}V_{(1)} \cdot  p({\bf U}_{(2)};{\bf V}_{(2)};w)
$$
where $(U,V,w)$ is the associator. Here Sweedler's notation is extended so as to mean
that the sum is taken over all partitions of the sequences {\bf U} and {\bf V} into pairs of
subsequences, ${\bf U} = {\bf U}_{(1)} \cup {\bf U}_{(2)}$ and 
${\bf V} = {\bf V}_{(1)} \cup {\bf V}_{(2)}$ such that 
$|{\bf U}_{(1)}|+|{\bf V}_{(1)} |  \geq 1,$ ${\bf U}_{(2)}\neq \emptyset ,$  ${\bf V}_{(2)}\neq \emptyset;$
 the expressions $U_{(1)}$ and $V_{(1)}$ are the right-normed products of
the elements of ${\bf U}_{(1)}$ and ${\bf V}_{(1)}$ respectively.

Our next step is to define $\tau$-analogues of these operations when $\tau^2=$ id.

Denote by
$$
[x,y]_{\tau}=(x\otimes y){\bf m}-(y\otimes x)\tau {\bf m}
$$
the {\em braided commutator}.
\begin{lm}\label{l6.1} If $\t$ is involutive, then the braided commutator $[x,y]_{\tau}$ is a  primitive operation.
\end{lm}
\Proof Let $y_1,y_2$ be arbitrary primitive elements of $B$. Suppose that $(y_1\otimes y_2)^{\tau}=\sum  z_{(2)}\otimes z_{(1)}$, where $z_{(1)},z_{(2)}\in B$. Then
$$
[y_1,y_2]_{\tau}=(y_1\otimes y_2){\bf m}-(y_1\otimes y_2)^{\tau}{\bf m}=y_1y_2-\sum  z_{(2)} z_{(1)},
$$
$$
([y_1,y_2]_{\tau})\Delta^b
=y_1y_2\underline{\otimes } 1+y_1\underline{\otimes } y_2+(y_1\underline{\otimes } y_2)^{\tau} +1\underline{\otimes } y_1y_2
$$
$$
-\sum (z_{(2)}z_{(1)}\underline{\otimes } 1+z_{(2)}\underline{\otimes } z_{(1)}+(z_{(2)}\underline{\otimes } z_{(1)})^{\tau} +1\underline{\otimes } z_{(2)}z_{(1)}).
$$
Notice that $\sum (z_{(2)}\underline{\otimes } z_{(1)})^{\tau}=((y_1\otimes y_2)^{\tau})^{\tau}=(y_1\otimes y_2)^{\tau^2}=y_1\otimes y_2$ since $\tau$ is involutive. Consequently,
$$
([y_1,y_2]_{\tau})\Delta^b=[y_1,y_2]_{\tau}\otimes 1 + 1\otimes [y_1,y_2]_{\tau};
$$
that is, the braided commutator is a primitive operation. $\Box$

We define the following system of operations
$$
P=P_{m,n} : B^{\otimes m}\otimes B^{\otimes n}\otimes B=B^{\otimes (m+n+1)}\rightarrow B, \ \ m,n\geq 1,
$$
by induction on $m+n$. Put
$$
P_{1,1}=({\bf m}_1-{\bf m}_2){\bf m}.
$$
Notice that $P_{1,1}$ is the associator, i.e., $P_{1,1}(x,y,z)=(x,y,z)$. If $m+n>1$, then a recursive formula for $P_{m,n}$ is given by
$$
P_{m,n}=(R_m\otimes R_n\otimes \mathrm{id})P_{1,1}-
\sum_{r+s>0}\sum_{\pi\in \mathrm{Sh}_m^r}
\sum_{\s\in \mathrm{Sh}_n^s} [\pi][\s]_{(m)}\nu_m^{r+1,m+s} 
(R_r\otimes R_s\otimes P_{m-r,n-s}){\bf m}_1 \ {\bf m},
$$
where $\sum_{r+s>0}$ means the summation over all $r,s$ such that $r\geq 0, s\geq 0, m-r\geq 0, n-s\geq 0$, and $r+s>0$. Moreover, in this case we can assume that $m-r\geq 1, n-s\geq 1$ because  $P_{i,j}=0$ if $i+j<2$. 

 It is not difficult to check that in the case of the ordinary flip  the operations $P_{m,n}$ coincide with the above 
 primitive operations $p_{m,n}.$  The main result of \cite{SU1} states that the set of all operations 
 $p_{m,n},$ $m,n\geq 1$, together with the commutator $[\cdot,\cdot]$ represents 
 a complete and an independent system of primitive operations in nonassociative algebras.
 
\begin{theor}\label{t6.1} If $\tau$ is involutive, then the operations $P_{m,n}$ 
for all $m,n\geq 1$ are primitive.
\end{theor}
\Proof Let $B$ be an arbitrary braided nonassociative Hopf algebra with braiding $\tau$. 
Let $V$ be the space of primitive elements of $B$. 
Denote by $\Delta^b$  the coproduct of $B$.   Notice that $P_{m,n}$ is a primitive operation if and only if the equality 
\begin{equation} \label{f6.1}
P_{m,n}\Delta^b=P_{m,n}\otimes 1+1\otimes P_{m,n}
\end{equation}
holds on $V^{\otimes (m+n+1)}$. Using (\ref{f5.3}), this equality can be rewritten 
in the form \begin{equation} \label{f7.9}
F_1 P_1+F_2 P_2+\ldots+F_s P_s=0,
\end{equation}
where $F_i$ are elements of the monoid algebra of the monoid generated by $\tau_1,\ldots,\tau_{m+n}$ for all $i$ and $P_1,P_2,\ldots,P_s$ are independent tensor products of parenthesis functions. Consequently, (\ref{f7.9}) holds on $V^{\otimes (m+n+1)}$ if and only if
\begin{equation} \label{f7.10}
F_1=F_2=\ldots=F_s=0
\end{equation}
holds on $V^{\otimes (m+n+1)}$. 

Notice that if $\tau$ is an ordinary flip, then $P_{m,n}$ exactly becomes the primitive operation $p_{m,n}$ described in \cite{SU1}. Moreover, every algebra can be considered as a braided algebra with the ordinary flip. This means that the equalities (\ref{f7.10}) hold for the ordinary flip. 
Taking into account the {\it principle} (see Remark \ref{r2.1}), this implies that (\ref{f7.10}) is true for any involutive braiding $\tau$. Consequently, $P_{m,n}$ is a primitive operation. 
$\Box$

\section{Sabinin  algebras in symmetric categories}

\hspace*{\parindent}

Let $\tau$ be an involutive braiding.
D. Gurevich \cite{Gurevich} introduced the concept of a ``Lie $\tau$-algebra", which is a wide 
generalization of the notions of Lie, super Lie, and color Lie algebras, as follows.

A braided algebra  $L$ with involutive braiding $\tau$ and 
multiplication ${\bf m}: L\otimes L\rightarrow L$ is called a {\em Lie $\tau$-algebra}
if
$$
{\bf m}+\tau {\bf m}=0, 
$$
$$
(\mathrm{id}+\tau_1\tau_2+\tau_2\tau_1)({\bf m}\otimes \mathrm{id}){\bf m}=0.
$$
These identities are braided analogues of the antisymmetric identity  and  the Jacobi identity.

Lie algebras first appeared as tangent algebras of Lie groups. 
In the case of simply connected Lie groups, Lie algebras determine the corresponding 
groups up to isomorphism. It was shown by P.O.Miheev and L.V.Sabinin  \cite{MS} 
that a simply connected local analytic loop is determined up to isomorphism 
by a more sophisticated analogue of a tangent algebra with a series of multilinear operations. 
These algebras are now called  {\em Sabinin} algebras \cite{MPS14,MP10,Perez07}.

We are going to give a definition of Sabinin $\tau$-algebras. First we define analogues 
of the relations (\ref{f3.1}). Let $\tau$ be a braiding on a vector space $R$.  
An $m$-ary multilinear operation $M : R^{\otimes m}\rightarrow R$, 
$m\geq 1$, is called a {\em braided} (or $\tau$-) algebraic operation on $M$ if
\begin{equation} \label{f7.1}
(M\otimes \mathrm{id})\tau=\tau_m\tau_{m-1}\ldots\tau_1(\mathrm{id}\otimes M),
(\mathrm{id}\otimes M)\tau=\tau_1\tau_2\ldots\tau_m(M\otimes \mathrm{id}).
\end{equation}

 A vector space $A$ with an involutive braiding $\tau$  is called 
 a {\em Sabinin $\tau$-algebra} if it is endowed with the multilinear braided operations
$$
S_{m,2}(x_1,\ldots ,x_m,y,z),\ m\geq 0, 
$$
$$
\Phi_{m,n} (x_1,\ldots ,x_m,y_1,\ldots ,y_n),\ m\geq 1, n\geq 2,
$$
that satisfy the identities:
\begin{equation} \label{f7.2}
S_{m,2}+\tau_{m+1}S_{m,2}=0,
\end{equation}
\begin{equation} \label{f7.3}
S_{m+2,2}-\tau_{r+1}S_{m+2,2}+\sum_{k=0}^{r} \sum_{\pi\in \mathrm{Sh}_r^k}[\pi]
(\mathrm{id}^{\otimes k}\otimes S_{r-k,2}\otimes \mathrm{id}^{\otimes (m-r+2)}) S_{m-r+k+1,2}=0
\end{equation}
for all $1\leq r<m+2$,
\begin{equation} \label{f7.4}
(1+\tau_{m+1}\tau_{m+2}+\tau_{m+2}\tau_{m+1})S_{m+1,2}
+\sum_{k=0}^{r} \sum_{\pi\in\mathrm{Sh}_r^k} [\pi](\mathrm{id}^{\otimes k}\otimes S_{m-k,2}\otimes \mathrm{id}) S_{k,2}=0,
\end{equation}
and
\begin{equation} \label{f7.5}
\Phi_{m,n}-\tau_r\Phi_{m,n}=0
\end{equation}
for all integers $r$ such that $1\leq r<m+n$ and $r\neq m$.

If $\tau =\theta $ is the ordinary flip, then every Sabinin $\tau$-algebra is 
a Sabinin algebra \cite[Section 1]{Perez07},  \cite[Section 4]{SU1}.

Let $V$ be an arbitrary (nonassociative) braided algebra over a field ${\bf k}$ of characteristic zero with involutive braiding $\tau$.
For any $x_1,x_2,\ldots ,x_m,y,z \in A$ we put
\begin{equation} \label{f7.6}
S_{0,2}(y,z)=-[y,z]_{\tau}
\end{equation}
and
\begin{equation} \label{f7.7}
S_{m,2}(x_1,x_2,\ldots ,x_m,y,z)=u(-P_{m,1}+\tau_{m+1}P_{m,1}),
\end{equation}
where $u=x_1\otimes x_2\otimes\ldots \otimes x_m\otimes y\otimes z,\ m\geq 1$.

Let $T_n$ be the monoid generated by $\tau_1,\tau_2,\ldots,\tau_{n-1}$. Since $\tau$ is involutive, we have a homomorphism from the symmetric group $S_n$ to $T_n$ such that $t_i\mapsto \tau_i$ for all $1\leq i<n$. We will denote the image of any $\pi\in S_n$ in $T_n$ by $\bar{\pi},$ $n\geq 1.$.

If $m\geq 1, n \geq 2$, then for any $x_1,x_2,\ldots ,x_m,y_1,y_2,\ldots ,y_n\in V$ we put
\begin{equation} \label{f7.8}
\Phi_{m,n}(x_1,x_2,\ldots ,x_m,y_1,y_2,\ldots ,y_n)
=\frac{1}{m!}\frac{1}{n!}\sum_{\pi \in S_m, \sigma \in S_n} u \bar{\pi}\bar{\sigma}_{(m)} P_{m,n-1},
\end{equation}
where $u=x_1\otimes\ldots\otimes x_m\otimes y_1\otimes\ldots\otimes y_n$.

\begin{theor}\label{t7.1}
Let $B$ be a braided algebra with involutive braiding $\tau $ over a field ${\bf k}$ of characteristic zero. 
The set  $B$
 is a braided Sabinin $\tau$-algebra with respect to the operations $($\ref{f7.6}$)$, 
 $($\ref{f7.7}$)$, and $($\ref{f7.8}$).$
\end{theor}
\Proof First we check that the operations (\ref{f7.6}), (\ref{f7.7}), and(\ref{f7.8}) are braided operations.
Let $M$ be one of these operations.  Suppose that $M$ is an $m$-ary operation.
 Using the relations (\ref{f3.1}), we can rewrite the relations (\ref{f7.1}) in the form (\ref{f7.9})
where $F_i$ are elements of the monoid algebra of the monoid generated by $\tau_1,\ldots,\tau_{m-1}$ 
for all $i$ and $P_1,P_2,\ldots,P_s$ are independent tensor products of parenthesis functions. 
Consequently, (\ref{f7.9}) holds if and only if (\ref{f7.10}) holds. 
It is easy to check that the relations (\ref{f7.1}) hold for any $m$-ary operation if $\tau$ is the ordinary flip. 
This means that (\ref{f7.9}) and (\ref{f7.10}) also hold for the ordinary flip. The {\it principle } 
(see Remark \ref{r2.1}) implies that (\ref{f7.10}) is true for an arbitrary involutive braiding.

Notice that each of the identities (\ref{f7.2})-(\ref{f7.4}) also can be written in the form
(\ref{f7.9}), where $P_1,P_2,\ldots,P_s$ are independent parenthesis functions. Consequently,
this identity is equivalent to (\ref{f7.10}). Moreover, all identities (\ref{f7.2})-(\ref{f7.4}) 
hold for the ordinary flip  \cite[Section 4]{SU1}. This means that (\ref{f7.10}) holds for the ordinary flip.
Remark \ref{r2.1} again implies that the identities (\ref{f7.2})-(\ref{f7.4}) hold for any involutive braiding. $\Box$

\begin{co}\label{c7.1}
Let $B$ be a braided (nonassociative) Hopf algebra with involutive braiding 
$\tau $ over a field ${\bf k}$ of characteristic zero. 
The set $\mathfrak A$ of all primitive elements  of $B$
is a braided Sabinin $\tau$-algebra with respect to the operations $($\ref{f7.6}$)$, 
$($\ref{f7.7}$)$, and $($\ref{f7.8}$).$
\end{co}
\Proof 
By the theorem above the algebra $B$ itself is a $\tau$-Sabinin algebra with respect to the operations
$S_{m,2}$ and $\Phi _{m,n}.$ By theorem \ref{t6.1} these operations are primitive, consequently 
the set $\mathfrak A$ is closed with respect to all Sabinin $\tau $-operations.
$\Box$

\

Finally, we formulate two interesting open problems.

1. Let ${\bf k}\{ V\} $ be a free nonassociative algebra over a braided space $V$
with involutive braiding. 
%
Does  $V$ generate  the space {\rm Prim}$({\bf k}\{ V\} )$ of all primitive nonassociative polynomials
as a $\tau$ -Sabinin algebra?


This problem has positive solution for nonassociative product with ordinary flip
(see \cite[Corollary 3.3]{SU1} or \cite[Theorem 8.2]{Khar15}) and also for the case
of associative product with arbitrary involutive brading (see \cite[Theorem 7.6]{Khar15}).

For each nonassociative braided algebra $B$ with an involutive braiding, $B^{(s)}$  
denotes the  $\tau $-Sabinin algebra $B$ with respect to the operations 
$($\ref{f7.6}$)$,  $($\ref{f7.7}$)$, and $($\ref{f7.8}$).$

2. Is it true that each $\tau $-Sabinin algebra can be isomorphically embedded into a $\tau $-Sabinin algebra 
$B^{(s)}$ for a suitable $\tau $-algebra $B$?

Recall that a similar problem has a positive solution for Sabinin algebras
(when $\tau $ is the ordinary flip, see \cite{Perez07}). 
The problem has an affirmative answer in the case of associative product
with arbitrary symmetry $\tau $ because each $\tau $-Lie algebra is a subagebra 
of its associative universal enveloping $\tau $-algebra (see \cite[Theorem 7.3]{Khar15}).

The authors are obliged to the referee for the careful reading and important corrections.

\end{document}